\documentclass[12pt]{article}

\usepackage{fullpage}
\usepackage{graphicx}
\usepackage{epstopdf}
\usepackage{float}
\usepackage[utf8]{inputenc}
\usepackage{t1enc}

\usepackage{amssymb,amsmath}
\usepackage[centerlast]{subfigure}
\usepackage{wrapfig}
\usepackage{amsthm}
\usepackage{color}
\usepackage{xcolor}
\usepackage{graphicx}
\usepackage{authblk}

\usepackage{hyperref}

\newtheorem{theorem}{Theorem}
\newtheorem{lemma}{Lemma}

\theoremstyle{definition}

\theoremstyle{remark}

\def\rev#1{\textcolor{black}{#1}}

\title{Qualitative properties of space-dependent SIR models with constant delay and their numerical solutions}
\author[1,2,3]{B.M. Takács\footnote{Corresponding author. E-mail: takacs.balint.mate@gmail.com}}
\author[1,2,3]{I. Faragó}
\author[2,3]{R. Horváth}
\author[4]{D. Repov\v s}

\affil[1]{Institute of Mathematics, Eötvös Loránd University, Budapest, Hungary}

\affil[2]{Institute of Mathematics, Budapest University of Technology and Economics, Budapest, Hungary}

\affil[3]{MTA-ELTE NumNet Research Group, Budapest, Hungary}

\affil[4]{Faculty of Education, Faculty of Mathematics and Physics at University of Ljubljana, and Institute of Mathematics, Physics and Mechanics, Ljubljana, Slovenia}

\date{ }

\begin{document}
\maketitle

\begin{abstract}
In this article a space-dependent epidemic model equipped with a constant latency period is examined. We construct a delay partial integro-differential equation and show that its solution possesses some biologically reasonable features. We propose some numerical schemes and show that by choosing the time step to be sufficiently small the schemes preserve the qualitative properties of the original continuous model. Finally, some numerical experiments are presented that confirm the aforementioned theoretical results.
\end{abstract}



\section{Introduction}
The increasing rate of globalization in recent decades led to an interconnected world in which diseases can spread faster than ever: this threat became reality in 2020 with the outbreak of the Covid-19 virus. Thus, the importance of adequate modeling of epidemics is apparent. 

One of the most frequently used tools in mathematics to model the spread of diseases is the SIR model \cite{Capasso,km27}. These models can be used to describe the spread of some feature among a group of individuals. We split our population into three categories: class $S$ contains those who do not have the property yet, class $I$ includes those who have the property and they can also transmit it to others, and class $R$ has those who had the property, but they do not have it any more. Note that these models can be used not only in epidemics but also in several other fields, see e.g.  \cite{Alsenafi,Bonnasse-Gahot,Rendine,Volkening} - however, in this article we are going to focus on epidemic models. The SIR models are generally systems of ordinary differential equations that do not take into the account the spatial location of the individuals of the above groups, although Kendall proposed a possible extension almost 60 years ago \cite{kend57,kend65}. This extension results in a system of partial integro-differential equations that can be solved only numerically. The numerical solution can be produced by some appropriate methods but it is a natural requirement that it should possess the characteristic qualitative properties of the disease propagation process. These properties were investigated in the recent papers \cite{Csomos:SIR,fh17,TakacsHorvathFarago:SIR,TakacsHadjimichael:SIR}. The qualitative properties of diffusion SIR models, which are also used to model spatial disease propagation, can be found e.g. in \cite{Capasso, Ducasse, Wu}. 

Some diseases take some time to develop inside an infected individual, so they do not start to infect upon their infection but only after a short period of time (so-called latency period). For this, we generalize the previous models by using delay integro-differential equations. In this paper, we formulate the basic qualitative properties of these models, we show that the solution of the continuous problems possesses these features and give sufficient conditions that guarantee the properties to the numerical solutions.

In Section \ref{Sec:MathModels}, we construct our mathematical models, and then we prove some properties of the analytical solution in Section \ref{Sec:analsol}. Section \ref{Sec:spatial_disc} deals with the spatial discretization of the model, while Section \ref{Sec:timedisc} discusses the application of the explicit Euler and Runge-Kutta methods. Then some numerical results are shown in Section \ref{Sec:NumExp}.

\section{Construction of the mathematical model}\label{Sec:MathModels}
The first, and possibly best known mathematical formulation of the SIR models comes from Kermack and McKendrick \cite{km27}, in which the authors propose the following model:
\begin{equation}\label{eq:Kermack}
\left\{
\begin{aligned}
S'(t) & = - a S(t)I(t),  \\
I'(t) & =   a S(t)  I(t) - b I(t), \\
R'(t) & = b I(t), 
\end{aligned}
\right.
\end{equation} 
where the functions $S(t), I(t), R(t): [0,\mathcal{T}] \rightarrow\mathbb{R}$ ($\mathcal{T} \in \mathbb{R}^+$) denote the number of individuals in class $S$ (susceptible), class $I$ (infected) and class $R$ (recovered) at time $t$, respectively, while $a, b \in \mathbb{R}^+$ are parameters describing the rate of infection and recovery, respectively. In this paper, we will consider the more realistic model
\begin{equation}\label{eq:delayc}
\left\{
\begin{aligned}
\frac{\partial S(t, x, y)}{\partial t} &= - S(t, x, y) F_I(t-\sigma,x,y) - c S(t,x,y),  \\
 \frac{\partial  I(t, x, y)}{\partial t}& = 
  S(t, x, y)F_I(t-\sigma,x,y)  - bI(t, x, y), \\
\frac{\partial R(t, x, y)}{\partial t} &= bI(t, x, y) + c S(t,x,y)
\end{aligned}
\right.
\end{equation} 
that also contains a delay term besides the generalizations of papers \cite{Csomos:SIR,TakacsHorvathFarago:SIR,TakacsHadjimichael:SIR}. In \eqref{eq:delayc}, $X:(0,\mathcal{T}]\times \Omega\to \mathbb{R}$, $(t,x,y)\mapsto X(t,x,y)$  ($X\in\{ S,I,R\}$) denotes the spatial density of the class $X$ at the point $(x,y)\in \Omega\subset \mathbb{R}^2$ ($\Omega$ is a connected, bounded open set) and at the time instant $t$. The parameter $b$ is the same as in \eqref{eq:Kermack}. The term $c S(t,x,y)$ ($c \in \mathbb{R}^+$) describes the effect of possible vaccination, which means that it is possible for some individuals to get from class $S$ to class $R$ without entering class $I$. The infection of a given individual is caused by the infected individuals in its surroundings, due to the assumed latency period at a previous time instant, and this effect depends on the spatial position of the individuals. Here we suppose that one can only be infected by others in a $\delta$-radius ($\delta \in \mathbb{R}^+$) neighborhood around itself and the effect of the infection is described by a given non-negative\rev{, continuous and} bounded weight function $W$. The constant delay $\sigma>0$ is the length of the latency period of the disease (see e.g.
\cite{Cooke:SIR_delay,MaSongTakeuchi:SIR_delay,XuMa:SIR_delay}). This phenomenon is included into the model by the delay term 
\begin{equation}\label{eq:F_I_def}
F_I(t-\sigma,x,y) = \iint_{B_\delta(x,y)} W(x',y') I(t-\sigma,x',y')\, dx' dy',\end{equation}
where $B_\delta(x,y)$ is the $\delta$ radius open ball centered at $(x,y)$ and $W$ is the weight function. 
The natural births and deaths are not taken into consideration in the model.

Because \eqref{eq:delayc} is a delay system, to obtain a properly posed problem we also need the values of the function $I$ on the time interval $[-\sigma,0]$. This is the so-called history function that will be denoted by $I_h:[-\sigma,0]\times \Omega\to \mathbb{R}$, $(t,x,y)\mapsto I_h(t,x,y)$. Later we will show that some assumptions are needed on the history function to assure that the solution is continuous \rev{in space and time} and behaves as expected.
The history functions of $S$ ans $R$ are denoted similarly but these functions do not appear in the model.

The model does not have boundary condition in a classical sense but due to the integral in \eqref{eq:F_I_def} we assume that $I$ is equal to zero outside the domain $\Omega$.  

In the next section we show that our problem has a unique solution, which behaves in a biologically reasonable way.

\section{Properties of the analytic solution}\label{Sec:analsol}

\subsection{Existence, uniqueness and smoothness of the solution}

In this section we discuss the solvability of system \eqref{eq:delayc}. The key tool will be the method of steps introduced by Bellman \cite{Bellman:MoS}, which involves the splitting of our time interval $[0,\mathcal{T}]$ into smaller intervals with length $\sigma$ with the 
grid points $0,\sigma,2\sigma,3\sigma,\ldots,\mathcal{T}$. The $(k+1)$th element of this list will be denoted by $t_k$. Let us denote the solution on the time interval $(t_{k-1}, t_{k} ]$ by $S_{k}(t,x,y)$, $I_k(t,x,y)$ and $R_k(t,x,y)$. Then the delay differential equation \eqref{eq:delayc} on a given interval $(t_k, t_{k+1}]$ ($k\geq1$) becomes a classical partial integro-differential equation
\begin{equation}\label{eq:bellmanstep}
\left\{
\begin{aligned}
\frac{\partial S_{k+1}(t, x, y)}{\partial t} &= - S_{k+1}(t, x, y) F_{I_{k}}(t-\sigma,x,y) - c S_{k+1} (t,x,y),  \\
 \frac{\partial  I_{k+1}(t, x, y)}{\partial t}& = 
  S_{k+1}(t, x, y)F_{I_{k}}(t-\sigma,x,y)  - bI_{k+1}(t, x, y), \\
\frac{\partial R_{k+1}(t, x, y)}{\partial t} &= bI_{k+1}(t, x, y) + c S_{k+1} (t,x,y), 
\end{aligned}
\right.
\end{equation} 
with initial conditions
$$S_{k+1}(t_k, x, y) = S_k (t_k, x, y), \quad I_{k+1}(t_k, x, y) = I_k (t_k, x, y), \quad R_{k+1}(t_k, x, y) = R_k (t_k, x, y). $$
In the case of $t \in (0, \sigma]$, \eqref{eq:bellmanstep} has the form
\begin{equation}\label{eq:bellmanstep_first}
\left\{
\begin{aligned}
\frac{\partial S_1(t, x, y)}{\partial t} &= - S_1(t, x, y) F_{I_{h}}(t-\sigma,x,y) - c S_{1} (t,x,y),  \\
 \frac{\partial  I_1(t, x, y)}{\partial t}& = 
  S_1(t, x, y)F_{I_{h}}(t-\sigma,x,y)  - bI_1(t, x, y), \\
\frac{\partial R_1(t, x, y)}{\partial t} &= bI_1(t, x, y) + c S_{1} (t,x,y), 
\end{aligned}
\right.
\end{equation} 
with initial conditions
$$S_{1}(0, x, y) = S_h (0, x, y), \quad I_{1}(0, x, y) = I_h (0, x, y), \quad R_{1}(0, x, y) = R_h (0, x, y). $$
Because the structures of \eqref{eq:bellmanstep} and \eqref{eq:bellmanstep_first} are the same, we will only consider the case of \eqref{eq:bellmanstep} but all the statements also will be valid for \eqref{eq:bellmanstep_first}. 

If we proceed to solve the equation step by step (i.e. solving it on $(0,\sigma]$ then on $(\sigma, 2 \sigma]$ and so on) then the term $F_{I_{k}}(t-\sigma,x,y)$ is known, since it only depends on $I_k (t,x,y)$. It is also important to notice that the first equation of \eqref{eq:bellmanstep} does not depend on the later ones, thus can be solved independently of the others. After this $S_{k+1}(t,x,y)$ will be known, which means that the only unknown is $I_{k+1}(t,x,y)$ in the second equation, so it can also be solved directly. Finally, the last equation can also be solved by integration, resulting in the solution of system \eqref{eq:bellmanstep}. 

For the sake of simplicity, we introduce the function $\tilde{F}_{k+1}(t,x,y):=F_{I_{k}}(t-\sigma,x,y).$ In other words, we shift the domain of $F_{I_{k+1}}(t-\sigma,x,y)$ from $(t_{k-1},t_{k}] \times \Omega$ to $(t_k, t_{k+1}] \times \Omega$. Then \eqref{eq:bellmanstep} has the form ($t \in (t_k, t_{k+1}]$)
\begin{equation}\label{eq:bellmanstep_tilde}
\left\{
\begin{aligned}
\frac{\partial S_{k+1}(t, x, y)}{\partial t} &= - S_{k+1}(t, x, y) \tilde{F}_{k+1}(t,x,y) - c S_{k+1} (t,x,y),  \\
 \frac{\partial  I_{k+1}(t, x, y)}{\partial t}& = 
  S_{k+1}(t, x, y)\tilde{F}_{k+1}(t,x,y)  - bI_{k+1}(t, x, y), \\
\frac{\partial R_{k+1}(t, x, y)}{\partial t} &= bI_{k+1}(t, x, y) + c S_{k+1} (t,x,y). 
\end{aligned}
\right.
\end{equation} 

\begin{theorem}\label{th:existence}
Assume that the history functions $S_h(t,x,y)$, $I_h(t,x,y)$ and $R_h(t,x,y)$ are continuous in time \rev{and also in spatial variables $(x,y)$}. Then system \eqref{eq:delayc} has a unique solution, which is continuously differentiable in time on $(0, \mathcal{T}]$ \rev{and is continuous in space}.
\end{theorem}

\begin{proof}
We can express the solution of \eqref{eq:bellmanstep_tilde} as
\begin{equation}\label{eq:bellmansolution}
\left\{
\begin{aligned}
 S_{k+1}(t, x, y) &= \mathcal{K}_1 \exp \left( - c t - \int_{t_k}^{t} \tilde{F}_{k+1}(s,x,y) ds \right) ,  \\
  I_{k+1}(t, x, y)& = \mathcal{K}_2 e^{-b t} + e^{-bt} \int_{t_k}^t e^{bs} S_{k+1}(s,x,y) \tilde{F}_{k+1}(s,x,y) ds, \\
 R_{k+1}(t, x, y) &= \mathcal{K}_3 + b \int_{t_k}^t I_{k+1}(s,x,y) ds + c \int_{t_k}^t S_{k+1}(s,x,y) ds. 
\end{aligned}
\right.
\end{equation} 
For the solution to be continuous \rev{in time} at point $t_k$, we take 
$$\mathcal{K}_1=S_{k}(t_k,x,y) \; e^{c t_k}, \qquad \mathcal{K}_2=I_{k}(t_k,x,y) \; e^{b t_k}, \qquad \mathcal{K}_3=R_{k}(t_k,x,y),$$
since in this case $\lim_{t \rightarrow t_k} X_{k+1}(t,x,y) = X_{k}(t_k,x,y)$ for $X \in \{ S, I, R\}$. (Note that the values of $\mathcal{K}_1$, $\mathcal{K}_2$ and $\mathcal{K}_3$ are defined differently for different choices of $k$.) Also, because of the form of the solutions, it can be proved by induction that if the history functions are continuous \rev{in time}, then our solution will also be continuous \rev{in time} on the whole solution domain. \rev{It can be shown similarly that if the history functions are continuous in space, then the solution is also continuous in space since the weighting function $W$ was assumed to be continuous.}
We show that the solution is not only continuous \rev{in time} but \rev{is} also continuously differentiable \rev{in that variable}. By the fundamental theorem of calculus, it is also easy to see that the solutions are continuously differentiable in $t$ on each $(t_k,t_{k+1})$. At $t_k$ ($k \geq 1$), we show that 
$$\lim_{t \rightarrow t_k} \dfrac{\partial X_{k+1}(t,x,y)}{\partial t}  = \left. \dfrac{\partial X_{k}(t,x,y)}{\partial t}\right|_{t=t_k} \text{ for } X \in \{ S, I, R\}.$$
However, because of the previous arguments we know that the right-hand side of \eqref{eq:bellmanstep} is continuous \rev{in time} on $[t_k,t_{k+1})$, so e.g. in the case of $S_{k+1}$, we have
$$\lim_{t \rightarrow t_k} \dfrac{\partial S_{k+1}(t,x,y)}{\partial t} = \lim_{t \rightarrow t_k} \left( - S_{k+1}(t, x, y) \tilde{F}_{k+1}(t,x,y) - c S_{k+1} (t,x,y) \right)=$$ 
$$= - S_{k}(t_k, x, y) \tilde{F}_{k}(t_k,x,y) - c S_{k} (t_k,x,y) =\left. \dfrac{\partial S_{k}(t,x,y)}{\partial t}\right|_{t=t_k}.$$
A similar argument holds in the case of $I$ and $R$. This shows the continuous differentiability of the solution \rev{in variable $t$}.

Another question is whether the solution \eqref{eq:bellmansolution} is the only one. It can easily be shown that taking account only the first equation of \eqref{eq:bellmanstep_tilde}, it can be thought of as an ordinary differential equation, and its right-hand side fulfils the Lipschitz property. The same can also be shown for the second one (assuming that $S_{k+1}(t,x,y)$ is a known function), and also for the last one. Consequently, solution \eqref{eq:bellmansolution} is unique.
\end{proof}

In the next section we observe whether the continuous solution possesses some biological features.

\subsection{Qualitative properties of the analytic solution}
A mathematical model is considered to be reasonable not only when the system has only one solution but the behavior of such solution should also possess the properties of the biological model. Here we are going to observe four of these properties.

A natural requirement is that the density of each species cannot be negative.
\begin{itemize}
\item[$C_1$]: The functions $S(t,x,y)$, $I(t,x,y)$ and $R(t,x,y)$ should be non-negative.
\end{itemize}

Since we assume that there are no births or natural deaths in our region (or their rate is assumed to be equal) and the individuals do not move, the sum of the densities of the three classes should remain constant at each point, and this should also hold for the total number of individuals in the observed domain (which is given by the integral of the sum of the densities).

\begin{itemize}
\item[$C_2$]: The sum $S(t,x,y)+I(t,x,y)+R(t,x,y)$ should be constant in $t$ for all $(x,y)$ points inside our domain. Also, a consequence of this is that the following integral is also constant in time:
$$\int_{\Omega} S(t,x,y) + I(t,x,y) + R(t,x,y) \; dx \; dy.$$
\end{itemize}

Since there is no way one individual can get to class $S$, the density of it cannot increase.

\begin{itemize}
\item[$C_3$]: The function $S(t,x,y)$ cannot increase in $t$ at each point $(x,y)$.\end{itemize}

Also, no individual can get out of class $R$, so its density cannot decrease.

\begin{itemize}
\item[$C_4$]: The function $R(t,x,y)$ cannot decrease in $t$ at each point $(x,y)$.
\end{itemize}

Our goal is to show that the properties $C_1$ -- $C_4$ are satisfied by the solution of system \eqref{eq:delayc}.

\begin{theorem}
Suppose that properties $C_1$--$C_4$ hold for our history functions, which are also continuous \rev{in space as well as in time}. Then $C_1$--$C_4$ also hold for the solution of system \eqref{eq:delayc} on a time interval $(0,\mathcal{T}]$.
\end{theorem}

\begin{proof}
By Theorem \ref{th:existence}, system \eqref{eq:delayc} has a unique solution which is continuously differentiable \rev{in variable $t$} for $t \in (0, \mathcal{T}]$.

We will proceed by induction: first we prove the properties on $(t_0, t_1]$, then by supposing that they hold on $(t_{k-1}, t_{k}]$, we prove them on $(t_k, t_{k+1}]$. However, since the proof on the first interval does not differ that much from the proof on any arbitrary one, we present here only the proof of the latter one.

By adding up the equations of \eqref{eq:bellmanstep_tilde}, it is clear that property $C_2$ holds. 

For property $C_1$, let us consider system \eqref{eq:bellmanstep_tilde}. The solution of this system is \eqref{eq:bellmansolution}. From this form, it is evident that since $\mathcal{K}_1$ was chosen to be $S_{k}(t_k,x,y) e^{c t_k} \geq 0$, the values of the function $S_{k+1}(t,x,y)$ are also non-negative. (If $k=0$ then $S_1(0,x,y)=S_h(0,x,y)\geq 0$.) From the second equality of \eqref{eq:bellmansolution}, since $\mathcal{K}_2 \geq 0$, $S_{k+1}(s,x,y) \geq 0$ and $\tilde{F}_{k+1}(s,x,y) \geq 0$ for $s \in [t_k, t]$ (by the induction assumption), the relation $I_{k+1}(t,x,y) \geq 0$ is satisfied. The non-negativity of $R_{k+1}(t,x,y)$ can be easily seen from the third equality of \eqref{eq:bellmansolution}, since both $I_{k+1}(s,x,y)$ and $S_{k+1}(s,x,y)$ are non-negative for all $s \in [t_k, t]$, and $\mathcal{K}_3=R_k(t_k,x,y)$ is also non-negative (if $k=0$, then $R_h(0,x,y)$ is non-negative because $C_1$ holds for the history function by assumption). Thus, we proved that $C_1$ holds.

Property $C_3$ can be seen if we consider the form of the function $S_{k+1}(t,x,y)$ in \eqref{eq:bellmansolution}: it is evident that since $c>0$ and $\tilde{F}_{k+1}(s,x,y) \geq 0$, the function is decreasing (here we also use that $\mathcal{K}_1 \geq 0$, as mentioned before).

Property $C_4$ is a consequence of the non-negativity of $I_{k+1}(t,x,y)$ and $S_{k+1}(t,x,y)$, since the right-hand side of the third equation of system \eqref{eq:delayc} is non-negative.
\end{proof}

In the next sections we examine the semi-discretized and the fully discretized versions of \eqref{eq:delayc} and check whether their solutions satisfy the analogous versions of $C_1$ -- $C_4$.

\section{The spatially discretized models and their properties}\label{Sec:spatial_disc}
It is evident that Bellman's method can only be used in practice, when the integral of the history function $I_h$ at the time $t=-\sigma$ is known, which is usually not the case. Because of this, in the following sections we propose a numerical approach to approximate these analytic solutions.

Upon looking at system \eqref{eq:delayc}, it is evident that the most problematic part is the fact that it contains integrals on its right-hand side. In Section \ref{Sec:numint} we approximate the integral term $F_I(t-\sigma,x,y)$ by using some cubatures, and then in Section \ref{Sec:spatgrid} we define a spatial grid on the domain, approximating the partial differential equation by a class of ordinary differential equations, by defining a separate equation for each grid point.

\subsection{The properties of the cubature formula model}\label{Sec:numint}
We define some two-dimensional cubature formula on the disc $B_{\delta}(x,y)$ with positive weights to approximate the integral $F_I(t,x,y)$. Let us introduce a set of cubature points $(x+\eta_i, y+\xi_i)$ on the disc $B_{\delta}(x,y)$ and the positive cubature weights $w_i>0$ ($i=1, \ldots, p$). Here the values $\eta_i$, $\xi_i$ and $w_i$ might be also dependent on $(x,y)$ in the most general setting but in this article for simplicity we use the same cubature formula
\begin{equation}\label{eq:kozelites}
\begin{aligned}
\mathcal{F}_I(t-\sigma,x,y) =  \sum_{i=1}^p w_{i} W(x+\eta_i, y+\xi_i) I(t-\sigma,x + \eta_i,y + \xi_i) \approx F_I(t-\sigma,x,y)
\end{aligned}
\end{equation}
for every point of the domain. 

Using approximation \eqref{eq:kozelites} we obtain the system of partial differential equations
\begin{equation}\label{eq:trapeq}
\left\{
\begin{aligned}
\frac{\partial S(t, x, y)}{\partial t} &= -S(t, x, y)\mathcal{F}_I(t-\sigma,x,y) - c S(t, x, y), \\
\frac{\partial I(t, x, y)}{\partial t} &= S(t, x, y)\mathcal{F}_I(t-\sigma,x,y) - bI(t, x, y), \\
\frac{\partial R(t, x, y)}{\partial t} &= bI(t, x, y) + c S(t,x,y).
\end{aligned}
\right.
\end{equation}
It is clear that the arguments detailed in Section \ref{Sec:analsol} can be used similarly, which results in the following theorem.

\begin{theorem}
Assume that the history functions are continuous \rev{in space and time} and properties $C_1$ -- $C_4$ hold for them. Then system \eqref{eq:trapeq} has a unique solution, which is continuously differentiable \rev{in time, continuous in space} and also has properties $C_1$ -- $C_4$.
\end{theorem}

A natural question is the choice of the numerical approximation of the integral. In \cite{TakacsHadjimichael:SIR} two separate choices of cubatures were investigated. One of them, the Elhay-Kautsky cubature results in a uniform distribution of points on the unit disc, while the other, the Gauss-Legendre cubature (which involves a transformation of the integral to a unit square) results in a non-uniform distribution. Numerical experiments show that while the first one works well for polynomials, the second one is better for arbitrary nonlinear functions. Since we cannot guarantee that our function $I(t,x,y)$ is a polynomial, we are going to use here the latter one. For further details of the different methods, see \cite{TakacsHadjimichael:SIR}.

\subsection{The properties of the semi-discretized model}\label{Sec:spatgrid}

To obtain the numerical solution, we assume that our domain $\Omega$ is a rectangle with one vertex at the origin, i.e. $\Omega = (0,A) \times (0,B)$, $A, B \in \mathbb{R}^+$. Note that the following arguments also hold for domains in a more general form, but involve much more careful choice of spatial grids.

Let us discretize our rectangle shaped domain using the spatial grid
$$\mathcal{G} = \left. \left\lbrace  (x_k,y_l)  \in \Omega\, \right| \; x_k = (k-1) h_x, y_l = (l-1) h_y, 1\leq k \leq K , 1 \leq l \leq L \right\rbrace ,$$
supposing that $(K-1)h_x =A$ and $(L-1) h_y = B$. This grid consists of $K L$ points with spatial step sizes $h_x$ and $h_y$, and we approximate the continuous solutions by a matrix containing the values at these grid points.  

After this semi-discretisation, we get the following set of equations:
\begin{equation}\label{eq:trapeq2}
\left\{
\begin{aligned}
\frac{dS_{k,l}(t)}{dt} &= -S_{k,l}(t)\hat{\mathcal{F}}_{k,l}(t-\sigma,x_k, y_l) - c S_{k,l}(t), \\
\frac{dI_{k,l}(t)}{dt} &= S_{k,l}(t)\hat{\mathcal{F}}_{k,l}(t-\sigma,x_k, y_l) - bI_{k,l}(t), \\
\frac{dR_{k,l}(t)}{dt} &= bI_{k,l}(t) + c S_{k,l}(t),
\end{aligned}
\right.
\end{equation}
in which $X_{k,l}(t)$ ($X \in \lbrace S,I,R\rbrace$) denotes the approximation of the function at grid point $(x_k,y_l)$ and at time $t$ and 
\begin{equation*}
\hat{\mathcal{F}}_{k,l}(t-\sigma,x_k, y_l) =
\sum_{i=1}^p w_{i} \; W(x_k+\eta_i, y_l+\xi_i) \hat{I}(t - \sigma,x_k + \eta_i,y_l + \xi_i).
\end{equation*}
Note that the points $(x_k + \eta_i, y_l + \xi_i)$ might not be part of $\mathcal{G}$, so there are no $I_{k,l}$ values assigned to them. Because of this, we will approximate them by some interpolation method using the nearest known $I_{k,l}$ values and positive coefficients. (This is the reason for the $\hat{I}$ notation.) In order to satisfy the qualitative properties, it is important to choose such interpolations that preserve non-negativity in the sense that if the known values are non-negative, then the function we get at the end of our process should also be non-negative. Such interpolations include monotone interpolation that uses piecewise cubic Hermite interpolating polynomials \cite{DoughertyEdelmanHyman:interpolation, FritschCarlson:interpolation} ('pchip' for short), which will be used in the numerical experiments.

As in the previous section, the methods described in Section \ref{Sec:analsol} can be used again for system \eqref{eq:trapeq2}, which results in the following theorem.

\begin{theorem}
Assume that the history functions are continuous  and properties $C_1$ -- $C_4$ hold for them. Then system \eqref{eq:trapeq2} has a unique  solution, which is continuously differentiable \rev{in time and} also has properties $C_1$ -- $C_4$.
\end{theorem}

In Section \ref{Sec:timedisc} we present two different numerical methods for system \eqref{eq:trapeq2}: first we solve it using the explicit Euler method via the Elsgolts approach \cite{els}, and later positivity-preserving Runge-Kutta methods \cite{Shu:TVD, Shu:TVD2, Shu:TVD3}.

\section{The fully discretized models and their properties}\label{Sec:timedisc}
\subsection{Application of the explicit Euler method}\label{Sec:elsg}
One of the key elements in the solution of delay differential equations is the fact that the discontinuities should be included in the mesh of the time discretisation. Since our history function is smooth, the only discontinuities in the higher order derivatives might appear at the points $k \sigma$, $k \in \mathbb{N}$. Based on this, we define them as
$$\mathcal{G}_t = \left\lbrace  t_{n/m}  \in [-\sigma, \mathcal{T}]\, \left| \; t_{n/m} = n \dfrac{\sigma}{m}, \; n \in \mathbb{Z}, \; - m\leq n \leq  m \dfrac{\mathcal{T}}{\sigma} ,\right\rbrace \right. ,$$
where $m$ is a positive integer.

On this above mesh, we can define the scheme
\begin{equation}\label{eq:elsg}
\left\{
\begin{aligned}
S^{n + 1} &= S^n - \tau S^{n} \circ T^{n-m} - c \tau S^n, \\
I^{n + 1} &= I^n + \tau S^{n} \circ T^{n-m} - b\tau I^n, \\
R^{n + 1} &= R^n + b\tau I^n + c \tau S^n,
\end{aligned}
\right.
\end{equation}
where $0\leq n \leq  m \dfrac{\mathcal{T}}{\sigma}$ and $\tau = \dfrac{\sigma}{m}$. The symbol $\circ$ is the element-by-element or Hadamard product of the matrices and the $(k,l)$ element of the matrix $X^n$ ($X \in \{ S, I, R\}$) is the approximation of the value $X_{k,l}(t_{n/m})$, moreover the $(k,l)$ element of $T^{n-m}$ gives the approximation of $\hat{\mathcal{F}}(t_{n/m}-\sigma,x_k,y_l)$ by interpolating the elements of $I^{n-m}$.

Instead of analyzing the previous numerical method in terms of its convergence, we are going to observe how well the model describes the real-life processes: more precisely, whether our model preserves the discrete versions of qualitative properties $C_1$ -- $C_4$. From now on, we denote these by $D_1$ -- $D_4$.

Now we prove that for a sufficiently small time-step (or in other words, a sufficiently large $m$) properties $D_1$-$D_4$ hold for $n > 0$. 
\begin{theorem}\label{th:elsg_theor}
Suppose that properties $D_1$--$D_4$ hold for the history functions discretized on the grids $\mathcal{G}$ and $\mathcal{G}_t$. Property $D_2$ holds without restrictions. Furthermore, if we assume that the time step satisfies
\begin{align}\label{eq:cond_elsg}
\tau= \dfrac{\sigma}{m} \leq \min \left\lbrace \dfrac{1}{\bar{T} + c}, \dfrac{1}{b} \right\rbrace,
\end{align}
in which 
\begin{align}\label{eq:tildeT}
\bar{T} = \max_{(x_k,y_l) \in \mathcal{G}} M \sum_{i=1}^p w_{i} W(x_k+\eta_i, y_l + \xi_i),
\end{align}
and 
\begin{align*}
	M=\max_{(x_k,y_l) \in \mathcal{G}} \left\{S(0,x_k,y_l) + I(0,x_k,y_l) + R(0,x_k,y_l)\right\},
\end{align*}
then properties $D_1$, $D_3$ and $D_4$ also hold up to the step $n \leq m \dfrac{\mathcal{T}}{\sigma}$. 
\end{theorem}

\begin{proof}
The proof is similar to the one of \cite[Theorem~2]{TakacsHorvathFarago:SIR}.

We proceed through induction in $n$: first we prove that the properties hold for the first few steps, then we prove that if the properties are true up to some step $n$ then they also hold for step $n+1$. Since the proof for some arbitrary step and the first few ones are very similar, we only present here the one for the latter. Note that the aforementioned properties should be proved for every element in the matrices $S^{n+1}$, $I^{n+1}$ and $R^{n+1}$, but since they are all similar, here we present the proof for an arbitrary element.

By adding up all the equations of \eqref{eq:elsg} we get property $D_2$. The first equation in the scheme can also be rewritten as
$$ S_{k,l}^{n+1} = S_{k,l}^n (1 - \tau T_{k,l}^{n-m} - c \tau).$$
Properties $C_1$ and $C_3$ hold iff
$$0 \leq 1 - \tau T_{k,l}^{n-m} - c \tau \leq 1.$$
The higher bound is true if $T_{k,l}^{n-m}\geq 0$ but this holds because of either the assumption of the induction condition or the first assumption of the theorem. Also, the inequality of the lower bound can be rephrased as
$$\tau \leq \dfrac{1}{T_{k,l}^{n-m} + c}.$$
It is easy to see that because of the construction of $\bar{T}$ and properties $D_1$ and $D_2$ (which hold at $t_{(n-m)/m}$ by the induction condition), 
$$T_{k,l}^{n-m} \leq \bar{T} \qquad \forall k,l: (x_k,y_l) \in \mathcal{G},$$
which means that if \eqref{eq:cond_elsg} holds then property $D_1$ holds for $S^{n+1}$.

Also, since $S^n$, $I^n$ and $T^{n-m}$ all have non-negative elements, $D_1$ holds also for $I^{n+1}$ if \eqref{eq:cond_elsg} holds. By similar considerations, $D_1$ and $D_4$ also hold for $R^{n+1}$, which proves the statement.
\end{proof}

An important remark is that the step size must be in the form $\dfrac{\sigma}{m}$, which means that by condition \eqref{eq:cond_elsg}, the theoretically best step size is in the form $\dfrac{\sigma}{\tilde{m}}$, where
$$\tilde{m}=\min \left\lbrace m \in \mathbb{N}^+ \;\left|  \; \dfrac{\sigma}{m}< \min \left\lbrace \dfrac{1}{\bar{T} + c}, \dfrac{1}{b} \right\rbrace \right\rbrace.\right.$$

\subsection{Application of Runge-Kutta methods}
To achieve higher order convergence in our numerical approximation, we need to use higher order methods. One of the most widely used ones are the Runge-Kutta methods. Since in this paper we consider a constant delay, these methods are easily applicable.

In the sequel, the investigations will be based on the Shu-Osher form. Now we introduce the necessary notations used in the case of the classical ordinary differential equations, and then generalize them for the delay equations.

Consider a time-dependent problem in the form
\begin{equation}\label{eq:RK_eq}
u'(t)=F(u(t)), \qquad u(0)=u^0,
\end{equation}
and a Runge-Kutta method given in the Butcher form \cite{Butcher:Book} with coefficients $(a_{ij}) \in \mathbb{R}^{s \times s}$ and $\mathbf{b} \in \mathbb{R}^s$. Let $\mathcal{B}$ be the following matrix containing all the aforementioned coefficients of the method:
$$\mathcal{B} = \left( \begin{matrix}
(a_{ij}) & \mathbf{0} \\ 
\mathbf{b}^T & \mathbf{0} 
\end{matrix}  \right).$$
in which $\mathbf{0}$ is the zero vector with length $s$ (the number of the stages in the method). Let us also denote the $(s+1)$-dimensional identity matrix by $\mathbb{I}$. If there is such a number $r>0$ for which $(\mathbb{I} + r \mathcal{B})$ is invertible, then the explicit Runge-Kutta method for equation \eqref{eq:RK_eq} can be expressed in the canonical Shu-Osher form \cite{Shu:TVD, Shu:TVD2, Shu:TVD3}
\begin{equation}\label{eq:SO_rep}
\begin{aligned}
u_{(i)}^{n} & = v_i u^{n} + \sum_{j=1}^{i-1} \alpha_{ij} \left( u_{(j)}^{n} + \dfrac{\Delta t}{r} F\left( u_{(j)}^{n} \right) \right), \qquad i=1,\dots, s+1,\\
u^{n+1} & = u_{(s+1)}^{n},
\end{aligned}
\end{equation}
in which $\alpha_{ij}, v_i$ are real constants, $\Delta t$ is the time step and $u_n$ approximates $u(n \Delta t)$. Moreover, we have $\boldsymbol{\alpha}_r = (\alpha_{ij}) = r (\mathbb{I} + r \mathcal{B})^{-1} \mathcal{B} \in \mathbb{R}^{(s+1) \times (s+1)}$ and $\boldsymbol{v}_r =(v_i) = (\mathbb{I} + r \mathcal{B})^{-1} \bar{\boldsymbol{e}} \in \mathbb{R}^{s+1}$ (here $\bar{\boldsymbol{e}} = (1, 1, \dots , 1)^{T}$). The reason for using such representation is that on the right-hand side of \eqref{eq:SO_rep} we have the linear combinations of the steps of the explicit Euler method with time step $\Delta t /r$. This means that the qualitative properties of the explicit Euler method may be transmitted to the Runge-Kutta methods.

Depending on the values of the parameter $r$ we might have different representations of the same scheme. The Shu-Osher representation having the largest possible value of $r$ for which $(\mathbb{I} + r \mathcal{B})^{-1}$ exists and $\boldsymbol{\alpha}_r$ and $\boldsymbol{v}_r$ have non-negative components is called optimal, and we define 
$$\mathcal{C}:=\max \{ r \geq 0: \exists (\mathbb{I} + r \mathcal{B})^{-1} \text{ and } \boldsymbol{\alpha}_r\geq 0, \boldsymbol{v}_r  \geq 0\}.$$
$\mathcal{C}$ is called the \textit{SSP (strong stability preserving) coefficient} and this constant will be used in the scheme ($r = \mathcal{C}$). For further reading, see \cite{GottliebKetchensonShu:SSPbook}.

In our case, the problem \eqref{eq:trapeq2} can be written in the form
\begin{equation}\label{eq:RK_eq2}
u'(t)=F(u(t),u(t-\sigma)), \qquad u \text{ is given on } [-\sigma,0].
\end{equation}
Because of this, the explicit Runge-Kutta method applied to \eqref{eq:RK_eq2} takes the form
\begin{equation}\label{eq:RK}
\begin{aligned}
u_{(i)}^{n} & = v_i u^{n} + \sum_{j=1}^{i-1} \alpha_{ij} \left( u_{(j)}^{n} + \dfrac{\tau}{\mathcal{C}} F\left( u_{(j)}^{n}, u^{n-m} \right) \right), \qquad i=1, \dots, s+1,\\
u^{n+1} & = u_{(s+1)}^{n},
\end{aligned}
\end{equation}
($0 \leq n \leq m \mathcal{T} /\sigma$) with the notations of Section \ref{Sec:elsg}.  

The next theorem states that for a sufficiently small time step, properties $D_1$--$D_4$ hold for the above scheme.
\begin{theorem}\label{th:RK_cond}
Consider an explicit Runge-Kutta method in the form \eqref{eq:RK} with SSP coefficient $\mathcal{C}>0$ and applied to system \eqref{eq:trapeq2} with non-negative history function. Then $D_2$ holds, and the properties $D_1$, $D_3$ and $D_4$ also hold if the step-size satisfies
\begin{equation}\label{RK_cond}
\tau= \dfrac{\sigma}{m} \leq \mathcal{C} \min \left\{ \dfrac{1}{\bar{T} + c}, \dfrac{1}{b}\right\}.
\end{equation}
\end{theorem} 

\begin{proof}
We prove the theorem by induction in $n$ but since the proof for the first and an arbitrary step is similar, here we present only the latter one.
Also, we are proving the required properties for a fixed $k,l$ point, since if they hold for an arbitrary point then they also hold for all of them. To make notations simpler, in the next sections we are using the notation $\widetilde{X}^n=X_{k,l}^n$ for $X \in \{ S, I, R, T\}$ and $k,l$ are arbitrary fixed indices. 

By applying method \eqref{eq:RK} to our equation \eqref{eq:trapeq2}, we get the following scheme for an $i$th intermediate step ($1 \leq i \leq s+1$)
\begin{equation}\label{eq:RK_istage}
\begin{aligned}
			\widetilde{S}_{(i)}^{n} &= v_i \widetilde{S}^{n} + \sum_{j=1}^{i-1} \alpha_{ij}\left(\widetilde{S}_{(j)}^{n} -
			\frac{\tau}{\mathcal{C}}\left(\widetilde{S}_{(j)}^{n} \widetilde{T}^{n-m} +c \widetilde{S}_{(j)}^{n}\right)\right),\\
			\widetilde{I}_{(i)}^{n} &= v_i \widetilde{I}^{n} + \sum_{j=1}^{i-1} \alpha_{ij}\left(\widetilde{I}_{(j)}^{n} +
			\frac{\tau}{\mathcal{C}}\left(\widetilde{S}_{(j)}^{n} \widetilde{T}^{n-m} -b \widetilde{I}_{(j)}^{n}\right)\right),\\
			\widetilde{R}_{(i)}^{n} &= v_i \widetilde{R}^{n} + \sum_{j=1}^{i-1} \alpha_{ij}\left(\widetilde{R}_{(j)}^{n} +
			\frac{\tau}{\mathcal{C}}\left(b \widetilde{I}_{(j)}^{n} + c \widetilde{S}_{(j)}^{n}\right)\right). 
\end{aligned}
\end{equation}
Let $Z_{(i)}^{n}=\widetilde{S}_{(i)}^{n} + \widetilde{I}_{(i)}^{n} + \widetilde{R}_{(i)}^{n}$.
Then adding up equations \eqref{eq:RK_istage} yields
$$\boldsymbol{Z} = Z^{n} \boldsymbol{v}_{\mathcal{C}}  + \boldsymbol{\alpha}_{\mathcal{C}} \boldsymbol{Z},$$
that is
$$\boldsymbol{Z} = Z^{n} (\mathbb{I} - \boldsymbol{\alpha}_{\mathcal{C}})^{-1}   \boldsymbol{v}_{\mathcal{C}},$$
in which $\boldsymbol{Z}=(Z_{(1)}^{n},Z_{(2)}^{n},\dots, Z_{(s+1)}^{n})^{T}$ and $Z^{n}=\widetilde{S}^{n} + \widetilde{I}^{n}+ \widetilde{R}^{n}$. It is well known that for our Runge-Kutta method to be consistent, we need that $v_i + \sum_{j=1}^{s+1} \alpha_{ij} = 1$ for all $1 \leq i \leq s+1$ and consequently $(\mathbb{I} - \boldsymbol{\alpha}_{\mathcal{C}})^{-1}\; \boldsymbol{v}_{\mathcal{C}} = \bar{\boldsymbol{e}}$, which shows that $Z_{(s+1)}^n = Z^n$. Then because of $Z_{(s+1)}^{n}=Z^{n+1}$, we get that
$$Z^{n+1} = \widetilde{S}^{n+1} + \widetilde{I}^{n+1} + \widetilde{R}^{n+1} = \widetilde{S}^{n} + \widetilde{I}^{n} + \widetilde{R}^{n} = Z^n, \quad \forall n,$$
so property $D_2$ holds.

The remaining three properties require that $\widetilde{S}^{n+1}$, $\widetilde{I}^{n+1}$ and $\widetilde{R}^{n+1}$ are non-negative, and $\widetilde{S}^{n+1}$ is non-increasing while $\widetilde{R}^{n+1}$ is non-decreasing. Let us rewrite the first two equations in the internal stage \eqref{eq:RK_istage} as
\begin{equation}\label{eq:RK_istage2}
\begin{aligned}
			\widetilde{S}_{(i)}^{n} &= v_i \widetilde{S}^{n} + \sum_{j=1}^{i-1} \alpha_{ij}\widetilde{S}_{(j)}^{n} \left(1 -
			\frac{\tau}{\mathcal{C}}\left(\widetilde{T}^{n-m} +c \right)\right),\\
			\widetilde{I}_{(i)}^{n} &= v_i \widetilde{I}^{n} + \frac{\tau}{\mathcal{C}}\sum_{j=1}^{i-1} \alpha_{ij}\widetilde{S}_{(j)}^{n} \widetilde{T}^{n-m} +
			\left(1 - \frac{\tau}{\mathcal{C}}b\right)\sum_{j=1}^{i-1} \alpha_{ij} \widetilde{I}_{(j)}^{n}.
\end{aligned}
\end{equation}
By the definition of $\widetilde{T}^{n-m}$ and by the assumption of the induction, and since we are using non-negative cubature and positivity-preserving interpolation, $\widetilde{T}^{n-m}$ is also non-negative.

It is clear that $\widetilde{S}_{(1)}^{n}=v_1 \widetilde{S}^{n}$, $\widetilde{I}_{(1)}^{n}=v_1 \widetilde{I}^{n}$ and $\widetilde{R}_{(1)}^{n}=v_1 \widetilde{R}^{n}$ (by the Shu-Osher form of explicit Runge-Kutta methods) and they are non-negative. Also, from the form of \eqref{eq:RK_istage2} it is clear that the non-negative property of $\widetilde{S}_{(i)}^{n}$ and $\widetilde{I}_{(i)}^{n}$ (and also of $\widetilde{R}_{(i)}^{n}$) holds for every $2 \leq i \leq s+1$ if the condition 
\begin{equation}\label{RK_proofcond}
0 \leq 1 - \frac{\tau}{\mathcal{C}} b \quad \text{ and } \quad 0 \leq 1 -
			\frac{\tau}{\mathcal{C}}\left(\widetilde{T}^{n-m} +c \right) 
\end{equation}
holds. The first part of this condition is satisfied because of the second term of the right-hand side of \eqref{RK_cond}. For the second part of \eqref{RK_proofcond}, observe that by the definition of $\bar{T}$ in \eqref{eq:tildeT}, we know that
$$\widetilde{T}^{n-m} \leq \bar{T}. $$
Therefore the second part of condition \eqref{RK_proofcond} is also satisfied because of \eqref{RK_cond}, so $D_1$ holds. 

Moreover, since $\widetilde{T}^{n-m} \geq 0$, we get that
\begin{equation*}
 1 -\frac{\tau}{\mathcal{C}}\left(\widetilde{T}^{n-m} +c \right) \leq 1.
\end{equation*}
With this, we get the following estimate for $\widetilde{S}_{(i)}^{n}$ from the first equation of \eqref{eq:RK_istage2}
\begin{equation}\label{eq:RK_proof1}
			\widetilde{S}_{(i)}^{n} \leq v_i \widetilde{S}^{n} + \sum_{j=1}^{i-1} \alpha_{ij}\widetilde{S}_{(j)}^{n}.
\end{equation}
As mentioned before, to achieve consistency we need $v_i + \sum_{j=1}^{i-1} \alpha_{ij} = v_i + \sum_{j=1}^{s+1} \alpha_{ij} = 1$ for all $1 \leq i \leq s+1$. Let $1\leq p \leq s+1$ be the stage index for which the value of $\widetilde{S}_{(p)}^{n}$ is the largest, i.e. $\widetilde{S}_{(i)}^{n} \leq \widetilde{S}_{(p)}^{n}$ for all $1 \leq i \leq s+1$. Then taking $i=p$ inside \eqref{eq:RK_proof1} yields
\begin{equation*}
			\widetilde{S}_{(p)}^{n} \leq v_p \widetilde{S}^{n} + \sum_{j=1}^{p-1} \alpha_{pj}\widetilde{S}_{(j)}^{n} \leq v_p \widetilde{S}^{n} + \sum_{j=1}^{p-1} \alpha_{pj}\widetilde{S}_{(p)}^{n} \leq \left( 1- \sum_{j=1}^{p-1} \alpha_{pj}\right) \widetilde{S}^{n} + \sum_{j=1}^{p-1} \alpha_{pj}\widetilde{S}_{(p)}^{n}.
\end{equation*}
By rearranging, we have
\begin{equation*}
			\left( 1- \sum_{j=1}^{p-1} \alpha_{pj}\right)\widetilde{S}_{(p)}^{n} \leq \left( 1- \sum_{j=1}^{p-1} \alpha_{pj}\right)\widetilde{S}^{n}, 
\end{equation*}
that is
\begin{equation*}
    \widetilde{S}_{(p)}^{n} \leq \widetilde{S}^{n}.
\end{equation*}

Consequently, $\widetilde{S}_{(i)}^{n} \leq \widetilde{S}^{n}$ for all $1 \leq i \leq s+1$ and also $\widetilde{S}^{n+1} = \widetilde{S}_{(s+1)}^{n} \leq \widetilde{S}^{n}$ which gives property $D_3$.

By the third equation of \eqref{eq:RK_istage}, it is evident that $\widetilde{R}_{(i)}^{n} \geq \widetilde{R}^{n}$, hence $\widetilde{R}^{n+1}=\widetilde{R}_{(s+1)}^{n} \geq \widetilde{R}^{n}$, which gives property $D_4$.  
\end{proof}

We should also note here that like in the case of the explicit Euler method, we cannot use arbitrary values for time steps, but they should be in the form $\sigma / m$. Therefore, the theoretically best time step is $\sigma / \tilde{m}$, where
$$\tilde{m} = \min \left\lbrace m \in \mathbb{N}^+ \;\left|  \; \dfrac{\sigma}{m}< \mathcal{C} \min  \left\lbrace \dfrac{1}{\bar{T} + c}, \dfrac{1}{b} \right\rbrace \right\rbrace.\right.$$

\section{Numerical experiments}\label{Sec:NumExp}
In this section we present some numerical experiments to confirm our previous results. First we show that the bounds we got in Theorems \ref{th:elsg_theor} and \ref{th:RK_cond} are sharp in the sense that the use of bigger time steps results in qualitatively bad behavior. Then, in the second part we present some graphs on which the solutions are compared for different values of $\sigma$ and their qualitative properties are checked.

\subsection{Construction of the test problem}
\begin{figure}[!ht]
\centering{\includegraphics[height=4.4cm]{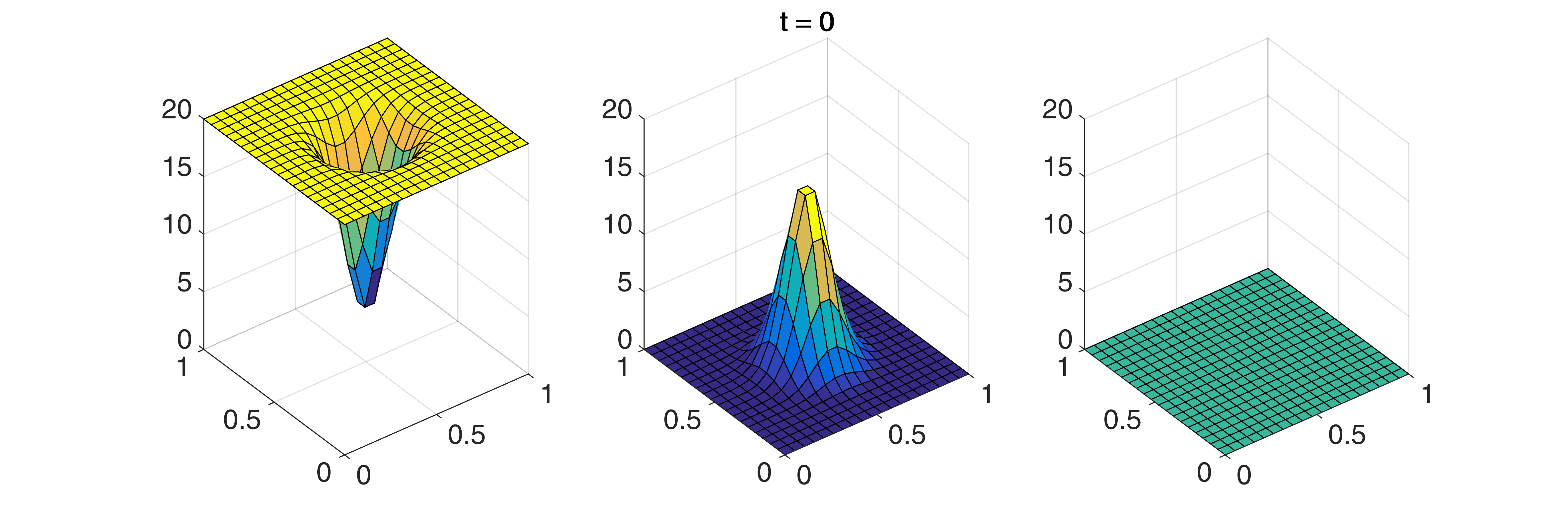}}
\parbox{14cm}{\caption{\label{fig:history}\footnotesize The history functions $S_h,I_h,R_h$ at $t=0$ shown in columns, respectively.}}
\end{figure}

The equation \eqref{eq:delayc} is solved on the rectangle domain $\Omega=(0,A)\times (0,B)$ ($A, B \in \mathbb{R}^+$). We set the parameters as $A=B=1$ and $c=0.01$. In order to be able to define the $F_I(t-\sigma, x,y)$ function, we choose the weight function
$$W(x',y') = a (-\| (x',y') - (x,y) \|_2+\delta).$$
In the numerical experiments, the choice $a=100$ is used. The radius $\delta$, delay-parameter $\sigma$ and the rate of recovery $b$ will be set later. The history functions are chosen as 
\begin{equation}\label{eq:initfunc}
\left\{
\begin{aligned}
 S_h(t, x, y) &= 20 - I_h(t,x,y),  \\
   I_h(t, x, y)&  
=  \frac{1}{2 \pi s^2}\exp\Bigg(-\frac{1}{2}\Bigg[\Bigg(\frac{x-\frac{1}{2} }{s}\Bigg)^2+\Bigg(\frac{y-\frac{1}{2} }{s}\Bigg)^2\Bigg]\Bigg) \left( 1 + \dfrac{t}{\sigma}\right), \\
R_h(t, x, y) &= 0,
\end{aligned}
\right.
\end{equation} 
for $t \in [-\sigma, 0]$, where $I_h(t,x,y)$ is a scaled Gaussian distribution with standard deviation $s=1/10$ concentrated at the middle of the domain $\left(1/2,1/2\right)$. Note that since $I_h$ is monotone increasing and continuous \rev{in time, continuous in space, is} non-negative and $1/(2 \pi s^2) \approx 15.92 < 20$, functions \eqref{eq:initfunc} fulfill properties $C_1$--$C_4$. The graphs of the history functions at time $t=0$ can be seen in Figure \ref{fig:history}.

The semi-discretization is carried out on a standard rectangular mesh with step sizes $h_x = h_y = 1/19$. As mentioned before, we can use different cubatures to approximate the integral $F_I(t,x,y)$ (see \cite{TakacsHadjimichael:SIR}) - we define the cubature as follows. First, we transform the disk-like infection domain with radius $\delta$ to the rectangle $[0,\delta]\times[0,2\pi)$ on the plane $(r,\vartheta)$ using a polar transformation with $x' = x + r \cos \vartheta$ and $y' = y + r \sin \vartheta$. Then, we transform this rectangle to the square $[0,1] \times [0,1)$ on the plane $(r',\theta')$ by using the linear transformation $r=\delta r'$ and $\vartheta=2\pi\vartheta'$ with the Jacobian determinant $2\pi\delta$. Using the transformations above, the integral $F_I(t-\sigma,x,y)$ has the form
\begin{equation*}
\int_0^1 \int_0^1 a (-r' \delta + \delta) I(t-\sigma, x + r' \delta \cos(2 \pi \vartheta') , y + r' \delta \sin(2 \pi \vartheta')) \; r' \; 2 \pi \delta^2 d r' d \vartheta'.
\end{equation*}
For the integration over the interior of the aforementioned square, we take the 40-point generalised Gaussian quadrature rule \cite{MaRokhlinWandzura:cubature} with the quadrature points by $\mu_1, \mu_2, \dots, \mu_{40}$ and corresponding weights $\omega_1 , \omega_2, \dots, \omega_{40}$. Therefore, the cubature has the form
\begin{equation*}
\begin{aligned}
\mathcal{F}_I (t-& \sigma, x, y) =\\
=\sum_{j=1}^{40} & \sum_{l=1}^{40} \omega_j \omega_l a (- \mu_j \delta + \delta) I\big(t-\sigma, x + \mu_j \delta \cos(2 \pi \mu_l), y + \mu_j \delta \sin(2 \pi \mu_l)\big) \mu_j 2 \pi \delta^2 = \\
=& \sum_{i=1}^{40^2} w_i W(x + \eta_i, y+\xi_i) I(t-\sigma, x+\eta_i, y+\xi_i)
\end{aligned}
\end{equation*}
where $i=40(j-1) + l$, $\eta_i = \mu_j \delta \cos(2 \pi \mu_l)$, $\xi_i = \mu_j \delta  \sin(2 \pi \mu_l)$ and $w_i = \omega_j \omega_l 2 \pi \delta^2 \mu_j$. Based on this, at the given spatial grid point $(x_k,y_l)$, the approximation  
\begin{equation*}
\hat{\mathcal{F}}_{k,l}(t-\sigma,x_k, y_l) =
\sum_{i=1}^p w_{i} \; W(x_k+\eta_i, y_l+\xi_i) \hat{I}(t - \sigma,x_k + \eta_i,y_l + \xi_i)
\end{equation*}
is used, where $\hat{I}(t - \sigma,x_k + \eta_i,y_l + \xi_i)$ is computed using piecewise cubic Hermite interpolation.

With the help of the previous constructions, the time discretization methods discussed in Section \ref{Sec:timedisc} now can be applied (see the next subsections).

\subsection{Sharpness of the time step bounds}
In the previous sections, namely in Theorems \ref{th:elsg_theor} and \ref{th:RK_cond} we gave sufficient conditions for the qualitatively good behavior of the numerical solution, i.e. if we use a smaller time step than the bounds, then our numerical solution possesses properties $D_1$ -- $D_4$. A natural question which might arise in the case of such sufficient conditions is the effect of the use of bigger time steps.

In Table \ref{Table1}, we can see the theoretical bound \eqref{eq:cond_elsg} 
\begin{equation}\label{eq:bound_comp}
\tau\le \min \left\{ \dfrac{1}{\left( 20 \sum_{j=1}^{40} \sum_{l=1}^{40} \omega_j \omega_l \; 100 \; \delta^3 (1 - \mu_j)  \mu_j 2 \pi \right) + 0.01} , 10\right\},
\end{equation}
denoted by "theor. b.", the actual time step in the form $\sigma / \tilde{m}$ (see the end of Section \ref{Sec:elsg}) denoted by "time step", and the bound calculated by experiments, i.e. the time step in the form $\sigma / m_{exp}$ for which the method works as expected but for $\sigma / (m_{exp}-1)$ it gives qualitatively inaccurate results - this value is denoted by "real b." in the table. Also, the difference $m_{exp} - \tilde{m}$ is denoted by "diff." in the table. The bound can be considered sharp when this value is zero. The last column shows the ratio of the "time step" and the "real bound", i.e. the sharpness of the bound we got from our theorem. As it can be seen, this ratio is $1$ for several parameter values. 

\begin{table}[ht]
\centering
\begin{tabular}{|l|l|l|l|l|l|l|}
\hline
$\delta$ & $\sigma$  & theor. b. & time step & real b. & diff. & ratio \\ \hline\hline
0.13 & 1 & 0.2169 & 0.2 & 0.2 & 0 & 1.0000  \\ \hline
0.12 & 1 & 0.2755 & 0.25 & 0.25 & 0 & 1.0000  \\ \hline
0.15 & 0.3 & 0.1413 & 0.1 & 0.1 & 0 & 1.0000  \\ \hline
0.15 & 0.5 & 0.1413 & 0.125 & 0.125 & 0 & 1.0000  \\ \hline
0.14 & 0.4 & 0.1737 & 0.1333 & 0.1333 & 0 & 1.0000  \\ \hline
0.13 & 0.5 & 0.2169 & 0.1667 & 0.1667 & 0 & 1.0000  \\ \hline
\end{tabular}

\caption{\label{Table1}
\footnotesize Numerical results for the explicit Euler method \eqref{eq:elsg} for various time steps with final time $\mathcal{T}=15$ and $b=0.05$ (the other parameters are given before).\bigskip}

\end{table}

\begin{figure}[!ht]
\centering{\includegraphics[height=6cm]{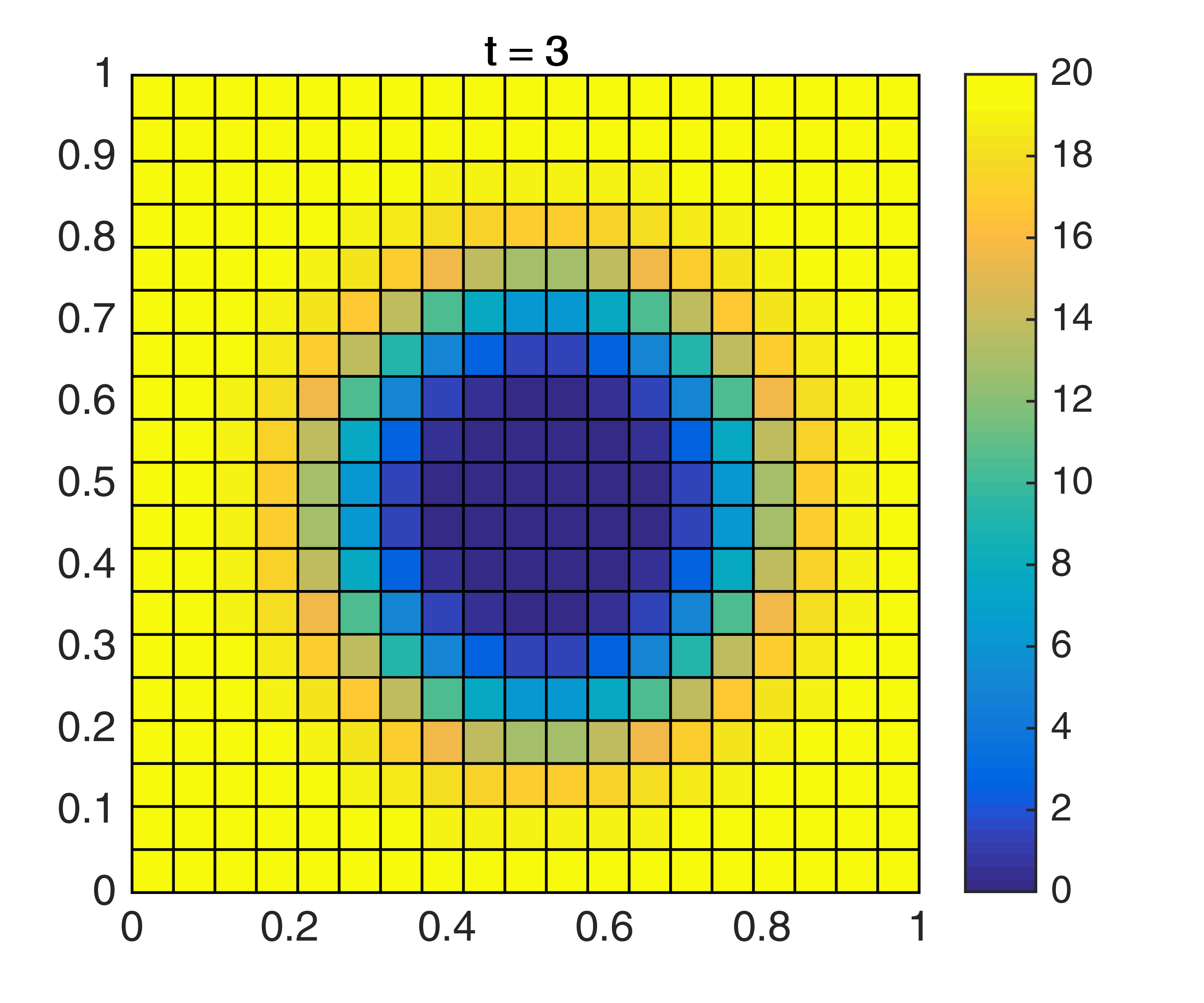}}
\centering{\includegraphics[height=6cm]{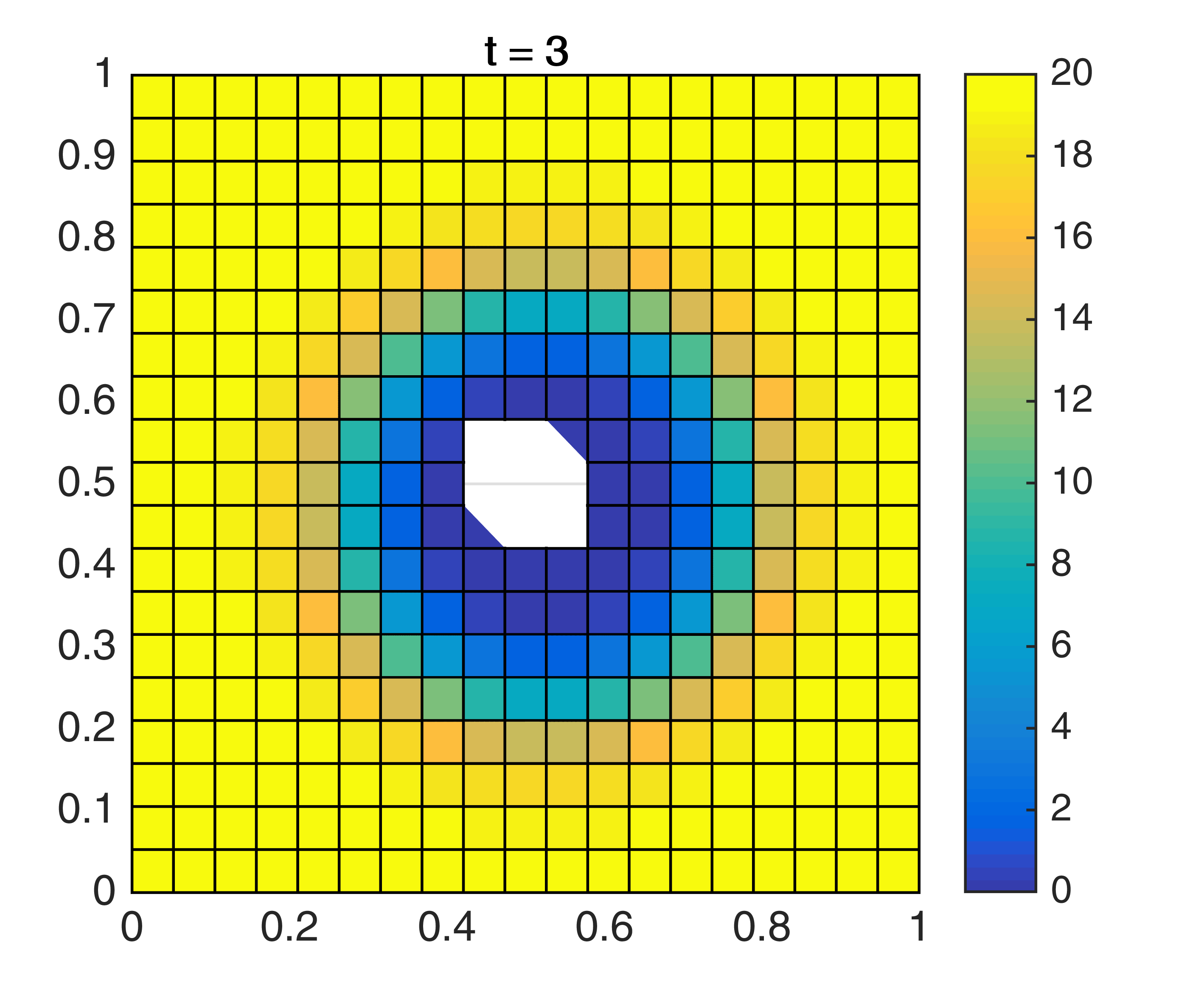}}
\parbox{14cm}{\caption{\label{fig:error}\footnotesize The numerical solution $S^n$ shown at final time $\mathcal{T}=3$ computed with time steps $\tau=0.25$ (left) and $\tau=0.333$ (right), with parameters $\delta=0.12$ and $\sigma=1$. The white area corresponds to those grid points at which the solution becomes negative.}}
\end{figure}

We show an example for the qualitatively bad behaviour of the method and show the sharpness of the obtained time step bound \eqref{eq:bound_comp} in the second row of Table \ref{Table1}. On the left panel of Figure \ref{fig:error}, the numerical solution $S^n$ can be seen at the time level $t=3$ obtained with time step $\tau=\sigma/4=1/4=0.25$. We can see that the solution is qualitatively correct, namely, the values are non-negative. But when we use the next possible time step $\tau=\sigma/3=1/3=0.3333$, we also get negative values. On the right panel of Figure \ref{fig:error} the white area corresponds to those grid points at which the solution is negative. Thus, the obtained bound is sharp.

The results for the second order Runge-Kutta method are presented in Table \ref{Table2}. Here we can see that the theoretical bound is not that sharp in all of the cases, so our condition is only sufficient, but not necessary. However, in some cases the theory gives time steps that are not far from the best possible one (when we have small numbers in the 'diff.' column). For different choices of the parameters and initial conditions, we might get even better results, resulting in a sharp bound even in the higher order case.

\begin{table}[ht]
\centering
\begin{tabular}{|l|l|l|l|l|l|l|l|}
\hline
$\delta$ & $\sigma$ & $b$ & theor. b. & time step & real b. & diff. & ratio \\ \hline\hline
0.13 & 1 & 0.1 & 0.2169 & 0.2000 & 0.5000 & 3 & 0.4000  \\ \hline
0.12 & 1 & 0.1 & 0.2755 & 0.2500 & 0.5000 & 2 & 0.5000  \\ \hline
0.13 & 0.5 & 0.05 & 0.2169 & 0.1667 & 0.2500 & 1 & 0.6667  \\ \hline
0.135 & 0.5 & 0.05 & 0.1937 & 0.1667 & 0.2500 & 1 & 0.6667  \\ \hline
0.135 & 0.4 & 0.01 & 0.1937 & 0.1333 & 0.2000 & 1 & 0.6667   \\  \hline
\end{tabular}

\caption{\label{Table2}
\footnotesize Numerical results for the RK2 method \eqref{eq:elsg} for various time steps with final time $\mathcal{T}=15$.\bigskip}

\end{table}

\subsection{Comparison of the cases with different delay parameters}
\begin{figure}[!ht]
\centering{\includegraphics[height=4.4cm]{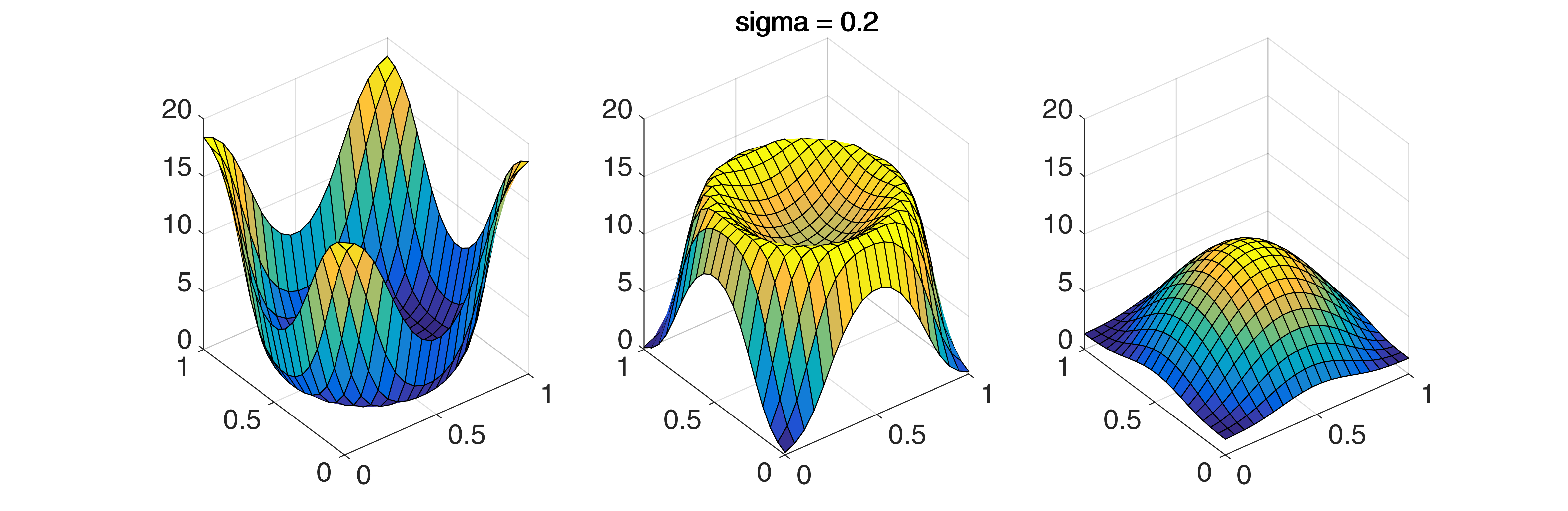}}
\centering{\includegraphics[height=4.4cm]{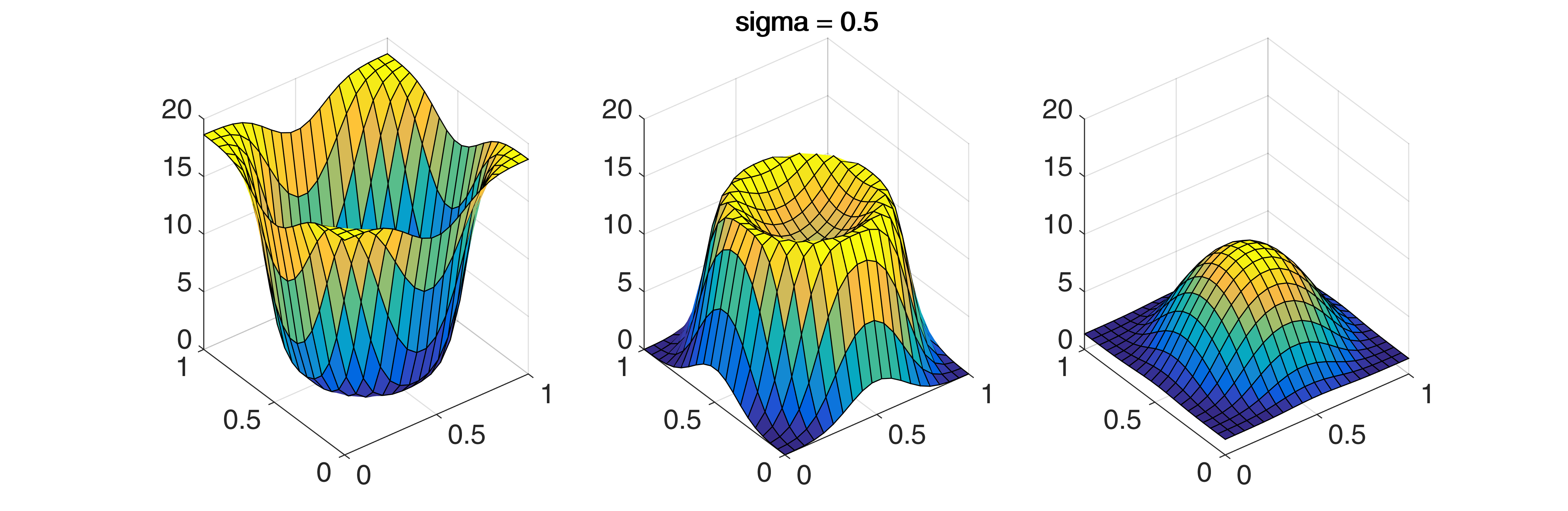}}
\centering{\includegraphics[height=4.4cm]{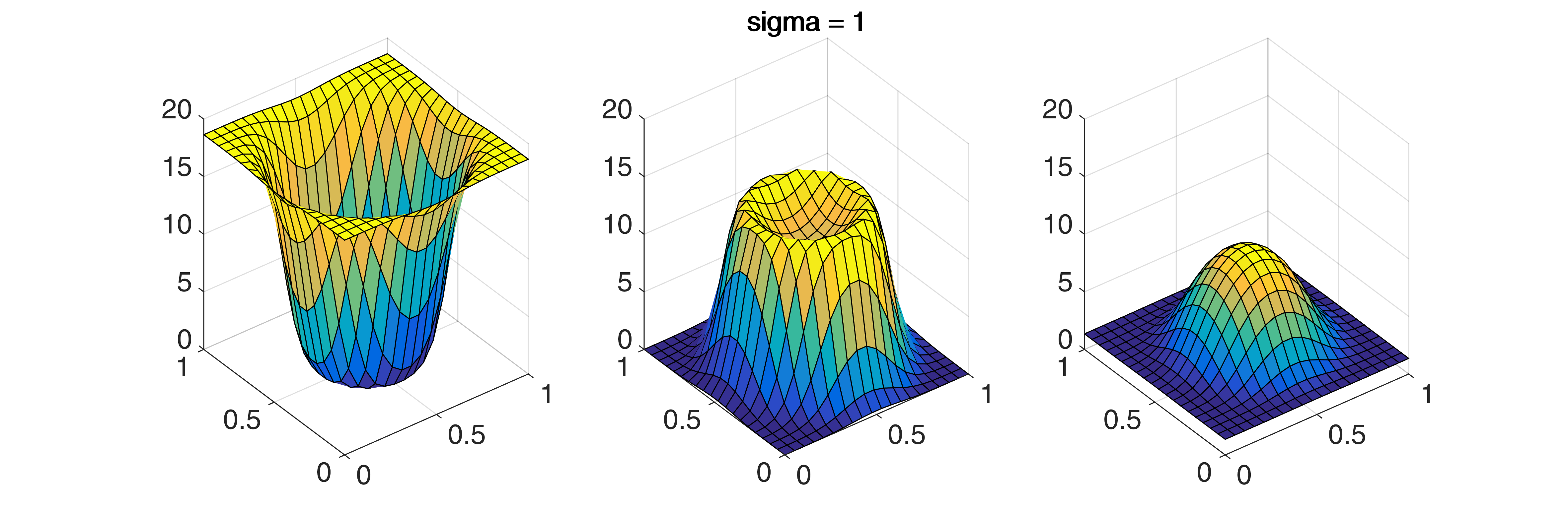}}
\centering{\includegraphics[height=4.4cm]{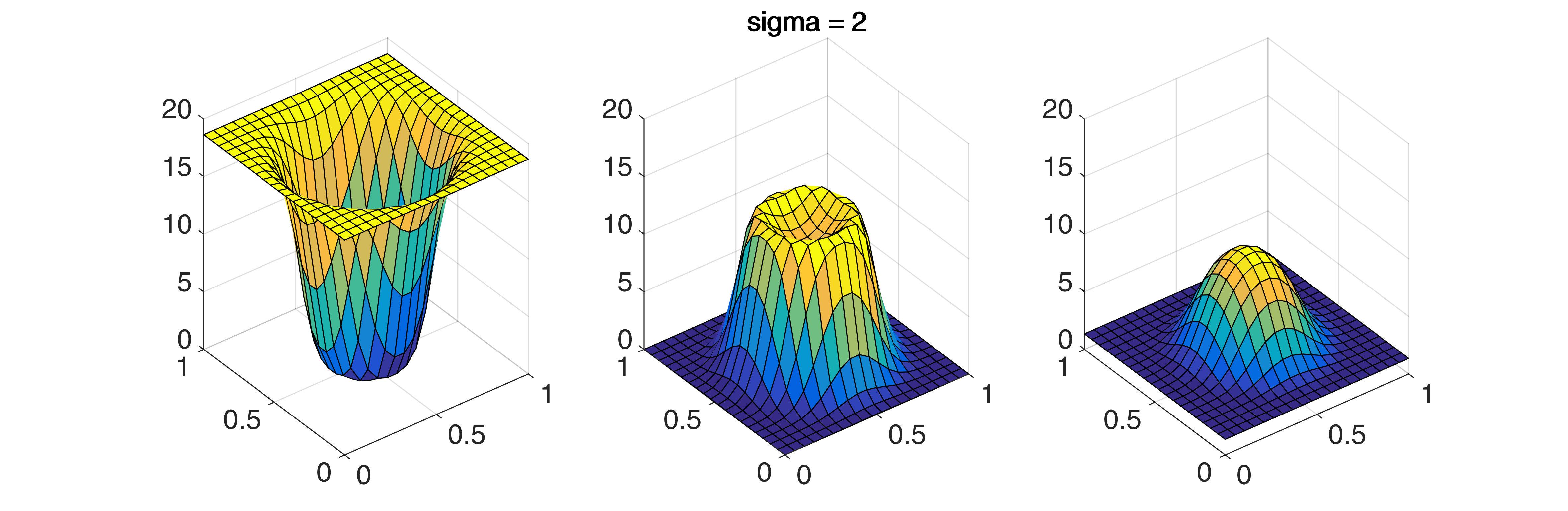}}
\parbox{14cm}{\caption{\label{fig:sir}\footnotesize The numerical solutions $S^n,I^n,R^n$ at $\mathcal{T}=7$ shown in columns, respectively, $\sigma$ values $\sigma=0.2$, $\sigma=0.5$, $\sigma=1.0$ and $\sigma=2.0$.}}
\end{figure}

In this section, we present some graphs of the numerical solutions of equation \eqref{eq:delayc} with different values of $\sigma$. We are plotting the solution at final time $\mathcal{T}=7$ with parameters $b=0.1$ and $\delta=0.1$, computed with the possible largest time step below the theoretical bound $0.4752$ (computed similarly as \eqref{eq:bound_comp} with $\mathcal{C}=1$) and the second order Runge-Kutta method is used. As we can see in Figure \ref{fig:sir}, the increase of parameter $\sigma$ results in a slower spread of the infection, which corresponds to the biological requirements (a longer latent period results in a slower pandemic).

\section{Conclusions, further work}
In this article we extended our previous works concerning space-dependent SIR models with the addition of constant delay to the equation. We could show using the method of steps that this new equation has a unique solution with biologically reasonable properties. Then, it was shown that numerical methods with carefully chosen step sizes will preserve the discrete versions of the aforementioned features.

One remaining question concerns the use of time steps not in the form $\sigma / m$, $(m \in \mathbb{N}^+)$. In these cases, an additional interpolation is needed to compute the values of our functions at times $t_{n/m} - \sigma$. However, if we use positivity-preserving interpolations (like the method 'pchip' mentioned in Section \ref{Sec:spatgrid}) then similar statements can be formulated as Theorems \ref{th:elsg_theor} or \ref{th:RK_cond}, although the computational time increases considerably.

Another possible extension of the previous methods include the introduction of non-constant delay, i.e. given by a time- or space-dependent function. In these cases the addition of an interpolation step to our algorithm is needed, and a more careful choice of time steps is also required.

\section*{Acknowledgements}
 The research by the authors B.M.T., I.F. and R.H. reported in this paper and carried out at BME has been supported by the NRDI Fund (TKP2020 NC, Grant No. BME-NC) based on the charter of bolster issued by the NRDI Office under the auspices of the Ministry for Innovation and Technology, and the Hungarian Ministry of Human Capacities OTKA
grant SNN125119.

The work of the author I.F. was completed in the ELTE Institutional Excellence Program (TKP2020-IKA-05) financed by the
Hungarian Ministry of Human Capacities.

The research of the author D.R. reported in this paper was supported by the Slovenian Research Agency grants P1-0292, N1-0114, N1-0083, N1-0064, and J1-8131.

\end{document}